\documentclass[preprint,12pt]{elsarticle}
\usepackage{graphicx}
\usepackage{natbib}
\usepackage{alltt}
\usepackage{amsmath}
\usepackage{amssymb}
\usepackage{algorithm}
\usepackage{algpseudocode}
\usepackage[hidelinks, pdftex]{hyperref}
\usepackage{float}
\makeatletter
\def\ps@pprintTitle{%
  \let\@oddhead\@empty
  \let\@evenhead\@empty
  \let\@oddfoot\@empty
  \let\@evenfoot\@oddfoot
}
\makeatother

\usepackage{natbib}
\usepackage{mathtools}
\usepackage{xcolor}

\begin{document}

\begin{frontmatter}
\title{Computationally Efficient Analysis of Energy Distribution Networks using Finite Volume Method and Interpolatory Model Order Reduction}
    \author[1]{Saleha Kiran \footnote{This work is funded by HEC NRPU No: 20-14782} }
\author[2]{Farhan Hussain}

\author[1]{Mian Ilyas Ahmad \corref{mycorrespondingauthor}}
\cortext[mycorrespondingauthor]{Corresponding author}
\ead{m.ilyas@sines.nust.edu.pk}

\address[1]{School of Interdisciplinary Engineering and Sciences, National University of Sciences and Technology Pakistan}
\address[2]{College of Electrical and Mechanical Engineering, National University of Sciences and Technology Pakistan}
\begin{abstract}
Energy distribution networks are crucial for human societies and since they often cover large geographical areas, their physical analysis is challenging. Modeling and simulation can be used to analyze such complex energy networks. In this paper, we performed discretization on the underlying partial differential equations of a pipeline and identified the complete network model for the gas, water and power distribution network. The discretized network model can be represented in a unified state-space form, that is, a specific matrix vector form for a set of differential-algebraic equations (DAEs). Due to the size and complexity of these models, their simulation can be computationally expensive. To address this issue, we applied model order reduction to the original unified mathematical model, constructing a reduced-order model that approximates the behavior of the original model with minimal computational cost. Specifically, we utilize the tangential iterative rational Krylov algorithm (tIRKA) which ensure that the interpolation condition of the linear part of the original and the reduced nonlinear unified models are interpolating along the predefined tangent directions. Specific scenarios of the energy networks are simulated with and without model order reduction and the behaviour is observed in terms of accuracy and computational cost. It is observed that with model order reduction, the computationally cost reduces significantly without compromising accuracy.


\end{abstract}

\end{frontmatter}

\section{Introduction}\label{s1}
We consider the problem of representing gas, water and power distribution networks in a unified mathematical framework as 
\begin{equation}\label{fom}
\begin{aligned}
    E\dot{x}(t)&=Ax(t)+Bu(t)+G(x,u)\\
    y(t)&=Cx(t)
    \end{aligned}
\end{equation}
where $E,A\in \mathbb{R}^{n\times n}$, $B\in \mathbb{R}^{n \times m}$, $C\in \mathbb{R}^{p \times n}$ are system matrices. Also $x(t)\in \mathbb{R}^n$, $u(t)\in \mathbb{R}^m$ and $y(t)\in \mathbb{R}^p$ are the state, input and output vectors, respectively. The function $G(x,u)$ is a nonlinear vector function of the state and input vectors. The matrix $E$ is singular and therefore the behaviour of these models for a specific input $u(t)$ require solution of the nonlinear differential and algebraic equations, which is computationally challenging for large value of $n$. A remedy to this problem is model order reduction, where the model is reduced to $r\ll n$ such that the reduced model is computationally cheap to simulate and has almost similar response as the full order model. Thus we are interested to compute a reduced system of the form
\begin{equation}\label{rom}
\begin{aligned}
    E_r\dot{x}_r(t)&=A_rx_r(t)+B_ru(t)+G_r(x_r,u)\\
    \tilde{y}(t)&=C_rx_r(t)
    \end{aligned}
\end{equation}
where $E_r,A_r\in \mathbb{R}^{r\times r}$, $B_r\in \mathbb{R}^{r \times m}$, $C_r\in \mathbb{R}^{p \times r}$ are the reduced system matrices. Also, $x_r(t)\in \mathbb{R}^r$ is reduced but $u(t)\in \mathbb{R}^m$ is the same and $\tilde{y}(t)\in \mathbb{R}^p$ is approximated but of the same size. The function $G_r(x_r,u)$ is the reduced nonlinear vector function of the reduced state and input vectors. 

The problem of model order reduction has been considered in detail for linear systems \cite{morAnt05}. Our focus is this work is its extension to nonlinear descriptor systems. In particular, we are utilizing the projection-based interpolation framework for the reduction of quadratic-bilinear descriptor systems as they are easily scalable to large-scale problems. In case of linear systems, the standard Arnoldi and Lanczos methods have been proposed it the literature, where the reduced method matches the Markov parameters of the original system. These frameworks have been extended to rational interpolation, where the transfer functions of the original and reduced system are implicitly matched at some predefined interpolation points by carefully constructing the basis matrices for r-dimensional subspaces. To obtain a good choice of interpolation points, the iterative rational Krylov algorithm has been introduced in \cite{morGugAB08}, that iteratively converges to fixed points where they are the mirror images of the eigenvalues of the reduced system.  The IRKA method on convergence satisfies the necessary conditions for optimality of the reduced system in terms of $\mathcal{H}_2$ norm. In case of multi-input multi-output systems, the concept of tangential interpolation has been utilized to compute the reduced system and its iterative version similar to IRKA has been proposed in the literature \cite{morGalVV04}. The issue, however, is that the tangential IRKA always constructs a strictly proper reduced system even though the original system is a descriptor system with both strictly proper and polynomial parts. This problem has been studied in \cite{morGugSW13}, where the descriptor system is first decomposed into strictly proper and polynomial parts and then reduction algorithm is applied only on the strictly proper part, retaining the polynomial part without reduction.  

The IRKA-type projection framework has been extended to bilinear, quadratic-bilinear and polynomial nonlinear systems where again basis matrices $V$ and $W$ are identified that implicitly match a series of frequency transfer functions associated with the nonlinear DAE system. Although these methods construct high-fidelity reduced models, there are many open questions that require further research. In particular, we observe how the nonlinear DAE system can be reduced without having undefined error at high frequencies. This is achieved by estimating the polynomial part of the nonlinear DAE system and retaining it in the feedforward matrix $D_r$ of the reduced system. The nonzero $D_r$ does not affect the tangential interpolation conditions by carefully constructing the state matrix $A_r$ and the input matrix $B_r$.   

The remaining part of the paper has been organized as follows. In Section~2, we present the mathematical modeling strategies for gas, water, and power distribution networks representing them in the unified system of DAEs \eqref{fom}. Section~3 discusses the proposed model order reduction framework for the DAE system. The numerical results are given in Section~4 and finally the conclusions and future directions.

\section{Energy Network Modelling}\label{s2}
In this section, we present the mathematical representation or modeling of the gas, water, and power distribution networks. Each of these networks are individually modeled and are expressed in the  unified mathematical framework \eqref{fom} for further analysis and design.
\subsection{Gas Distribution Network}
The gas dynamics in a single gas pipeline can be modeled as a one-dimensional isothermal Euler equation \cite{qiu2018efficient}, which for the case of constant temperature with spatial domain $[0,L]$ can be written as
\begin{equation}\label{1D_Euler_updated1}
\begin{aligned}
    \partial_t p &= - \frac{\gamma_0}{a} \partial_x q,\\
    \partial_t q &= -a\partial_x p-\frac{\gamma_0}{a}\partial_x \frac{q^2}{p}-ga\rho \partial_x h-\frac{\gamma_0\lambda(\varphi)}{2da}\frac{q|q|}{p},
    \end{aligned}
\end{equation}
where, $\partial_t$ denotes the partial derivative with respect to $t$ and $\rho$ represents the density of the gas. In addition, $p$ is the pressure of the gas and $\varphi$ is the volumetric flow rate written as $\varphi=\rho v$ in which $v$ is the velocity of the gas. Similarly, $d$ is the diameter, $g$ is the gravitational acceleration, $\gamma(T)$ is the specific heat ratio, and $h(x)$ is the elevation of the pipeline. In addition, $\lambda$ is the friction factor dependent on the flow rate of the gas, $T$ is the temperature of the gas and $z$ is the compressibility factor. The term $\frac{\gamma_0}{a}\frac{\partial}{\partial x} \frac{q^2}{p}$ is very small compared to $\partial_x p$ \cite{grundel2016numerical} and therefore can be neglected. In addition, the gravity term ${g\rho \frac{\partial h}{\partial x} }=\frac{gp \partial h}{\gamma_0 \partial x}$ can be ignored by assuming that pipes are buried underground at a homogeneous level. With these assumptions and $\gamma_0=c$,  the finite volume method is utilized as a 1D discretization problem, with boundary conditions $p=p_s,~\mbox{at}~x=0$ and $q=q_d~\mbox{at}~x=L$ to obtain
\begin{equation*}\label{final 1}
\begin{aligned}
     \frac{3h_{1}}{8}\frac{\partial q_{1}}{\partial t}+\frac{h_{1}}{8}\frac{\partial q_{2}}{\partial t}+\frac{a_1}{2}(p_2-p_1)=-\frac{c}{4}\frac{h_1\lambda_1}{a_1d_1}\frac{q_1|q_1|}{p_1}\\
    \frac{\partial p_i}{\partial t}(\frac{h_i+h_{i-1}}{2})+c (\frac{-1}{2a_{i-1}}q_{i-1}+(\frac{1}{2a_{i-1}}-\frac{1}{2a_i})q_i+\frac{1}{2a_i}q_{i+1})&=0\\
        \frac{(h_i+h_{i-1})}{2} \frac{\partial q_i}{\partial t}-\frac{a_{i-1}}{2}p_{i-1}+\frac{a_{i-1}-a_i}{2}p_i+\frac{a_i}{2}p_{i+1}=-\frac{c}{4}(\frac{h_{i-1}\lambda_{i-1}}{a_{i-1}d_{i-1}}&+\frac{h_{i}\lambda_{i}}{a_{i}d_{i}})\frac{q_i|q_i|}{p_i} \\ 
    \frac{h_{n-1}}{8}\frac{\partial p_{n-1}}{\partial t}+\frac{3h_{n-1}}{8}\frac{\partial p_{n}}{\partial t}+\frac{c}{2a_{n-1}}({q_d-q_{n-1}})=0
\end{aligned}
\end{equation*}
the Matrix-vector form can be written as
\begin{equation} \label{mat_vec_form}
        \begin{bmatrix} M_p &  \\ & M_q \end{bmatrix}	
        \begin{bmatrix}  \partial_t p \\\partial_t q \end{bmatrix}=
        \begin{bmatrix} 0 & K_{pq} \\ K_{qp}& 0 \end{bmatrix}
        \begin{bmatrix} p   \\ q \end{bmatrix}+
        \begin{bmatrix} B_q \\ 0 \end{bmatrix}q_d+
        \begin{bmatrix} 0  \\  B_p \end{bmatrix}p_s+
        \begin{bmatrix} 0  \\ g(p_s,p,q) \end{bmatrix}
\end{equation}
where by defining $h_{ij}=\frac{h_i+h_j}{2}$, the mass matrices $M_p$ and $M_q$ are
\begin{equation*}
\small
M_p\!=\!
\begin{bmatrix}
h_{12}& & &  \\
 &\!\!h_{23}& &  \\
 & &\!\!\ddots&  \\
 & & &\!\! h_{(n-2)(n-1)}\\
 & &&\!\!\!\!\frac{h_{n-1}}{8}&\!\!\!\!\frac{3h_{n-1}}{8}
 \end{bmatrix},~
M_q\!=\!
\begin{bmatrix}
\frac{3h_1}{8}&\frac{h_1}{8}& &  \\
 &\!\!h_{12}& &   \\
 & &\!\!h_{23}&  \\
 & & &\!\!\ddots& \\
 & & & &\!\! h_{(n-2)(n-1)}
 \end{bmatrix}.
\end{equation*}
Similarly by defining $\frac{1}{a_{ij}}=\frac{1}{a_i}-\frac{1}{a_j}$ and $\bar{a}_{ij}=a_i-a_j$, the coupling matrices $K_{pq}$ and $K_{qp}$ and the boundary matrices $B_p$ and $B_q$  in \eqref{mat_vec_form} are defined as
\begin{equation*}
\begin{aligned}
K_{pq}&=\frac{-c}{2}
\begin{bmatrix}
\frac{-1}{a_1}&\frac{1}{a_{12}}&\frac{1}{a_2}& &  \\
 &\ddots&\ddots&\ddots& \\
 & & \frac{-1}{a_{n-3}} &\frac{1}{a_{(n-3)(n-2)}}&\frac{1}{a_{n-2}}\\
 & & & \frac{-1}{a_{n-2}}&\frac{1}{a_{(n-2)(n-1)}}\\
 & & & &  \frac{-1}{a_{n-1}}
 \end{bmatrix}, \quad B_p= \frac{1}{2}\begin{bmatrix}a_1\\a_1\\0\\\vdots\\0\end{bmatrix}\!,\\
K_{qp}&=\frac{-1}{2}
\begin{bmatrix}
a_1& & &  \\
\bar{a}_{12}&a_2& & & \\
-a_2&\bar{a}_{23}&a_3& & \\
 & \ddots&\ddots&\ddots& \\
 &  &-a_{n-2}&\bar{a}_{(n-2)(n-1)}&a_{n-1}
 \end{bmatrix}, \quad B_q= \frac{1}{2}\begin{bmatrix} 0\\ \vdots\\0\\\frac{-c}{a_{n-1}}\\\frac{-c}{a_{n-1}}\end{bmatrix}\!\!.
 \end{aligned}
\end{equation*}
Finally the nonlinear vector function of frictional losses is written as 
\[
g(p_s,p,q)=-\frac{c}{4}\begin{bmatrix}(\frac{h_1\lambda_1}{a_1d_1}\frac{q_1|q_1|}{p_s})\\(\frac{h_1\lambda_1}{a_1d_1}+\frac{h_2\lambda_2}{a_2d_2})\frac{q_2|q_2|}{p_2}\\\vdots\\(\frac{h_{n-2}\lambda_{n-2}}{a_{n-2}d_{n-2}}+\frac{h_{n-1}\lambda_{n-1}}{a_{n-1}d_{n-1}})\frac{q_{n-1}|q_{n-1}|}{p_{n-1}}\end{bmatrix}.
\]
 In \eqref{mat_vec_form} the dynamics of a single pipeline are represented. In case of a network of pipelines, we can represent the system as a directed graph with edges $E$ and a set of nodes $N$. Edges are pipelines in a network while nodes are connections of pipelines. Nodes can be classified into three different types: supply node $N_s$, demand node $N_d$ and
interior nodes $N_i$. The supply nodes are the boundary nodes from which gas is
introduced in the network, while demand nodes are the boundary nodes
from which gas is extracted/withdrawn from the network. The interior node connects two or more edges of the network. The interior node that connects at least three
edges is called a junction node. One common simplification is to remove all interior nodes that are connecting exactly two edges and replace them with an edge. This
helps in reducing the number of algebraic constraints as only a single long edge will be part of the network until a junction node appears. For every $i$-th edge or pipeline that is part of the network's in flow, the pipeline has a supply pressure $p_s^{(i)}$ and a demand flow rate $q_d^{(i)}$ as defined by the boundary conditions of a pipeline. This means that at the junction nodes, the pressure of the outgoing pipe $p_s^{(j)}$  is equal to the pressure of the incoming pipe $p_s^{(i)}$, that is, $p_s^{(i)}=p_s^{(j)}$. In addition the mass flow rate of the incoming pipelines $ \sum_{k}q_d^{(i)}$ must be equal to the mass flow rate of the outgoing pipe $q_s^{(j)}$, that is, $\sum_{k}{q_d^{(i)}}=q_s^{(j)}$. In case of a two supply node and one demand node fork network having one junction node, the supply pressures are equal  $p_s^1=p_s^2=p_s^3$ and the demand mass flow rate at junction node is $q_d^1+q_d^2=q_d^3$.

Combining equation \eqref{mat_vec_form} and the constraints discussed earlier implies
\[
\begin{bmatrix}
M^{(1)}   &         &           & &    \\
          &M^{(2)}  &           & &    \\
          &         & M^{(3)}   & &    \\
          &         &           &0&    \\
          &         &           & & 0  \\
\end{bmatrix}
\frac{\partial }{\partial t}
\begin{bmatrix}
x^{(1)}  \\
x^{(2)}  \\
x^{(3)}  \\
q_d^{(4)}\\
q_d^{(5)}\\
\end{bmatrix}=
\begin{bmatrix}
K^{(1)}       &         &           &B_q^{(1)}&              \\
              &K^{(2)}  &           &         &B_q^{(2)}     \\
\beta_p^{(3)} &         & K^{(3)}   &         &              \\
              &         & e_3       &1        &1             \\
e_1           &e_2      &           &         &0            \\
\end{bmatrix}
\begin{bmatrix}
x^{(1)}  \\
x^{(2)}  \\
x^{(3)}  \\
q_d^{(4)}\\
q_d^{(5)}\\
\end{bmatrix}
\]\\

\[
+\begin{bmatrix}
B_p^{(1)}& \\
         &B_p^{(2)} \\
         & \\
         & \\
         &\\
\end{bmatrix}
\begin{bmatrix}
p_s^{(1)}  \\
p_s^{(2)}  \\
\end{bmatrix}
+\begin{bmatrix}
0  \\
0  \\
B_q^{(3)}  \\
0\\
0\\
\end{bmatrix}q_d^{(3)}
+\begin{bmatrix}
G_1{(x^{(1)},p_s^{(1)})}  \\
G_2{(x^{(2)},p_s^{(2)})}   \\
G_3{(x^{(3)},e_3^{(1)})}   \\
0\\
0\\
\end{bmatrix},
\]  
The state variables $x^{(i)}=[ p^{(i)}; q^{(i)}]$ represents pressure values and mass flow rates of the $i$-th pipe and $K^{(i)}$ are the coupling matrices for the $i$-th pipe. The last two rows in the state equations represent the algebraic constraints discussed above. The vector $\beta_p^{(3)}=[0,0,1,..,0,.. ]\otimes B_p^{(3)}$ incorporates equation \eqref{thirdpipe}, $e_1,~e_2,~e_3$ are elementary vectors with 1 or -1 at required places and 0's otherwise and the nonlinear vectors are $G_i{(x^{(i)},p_s^{(i)})}=[0; g(p_s^{(i)},p^{(i)},q^{(i)})]^T$. With these details, the required state space representation as in \eqref{fom} for gas distribution network can be written as
\begin{equation}\label{gas ssp}
    E_g\dot{x}_g(t)=A_gx_g(t)+B_gu_g(t)+G_g
\end{equation}

\subsection{Water Distribution Network}
Analogous to the modeling of gas distribution network, we first consider the representation of water flow through a single pipeline with spatial domain [0,L]. To simplify the mathematical representation, we utilize rigid water column theory \cite{izquierdo2004mathematical} which means that the following
continuity equation holds
 \begin{equation}
     \frac{\partial q}{\partial x}(x,t)=0.
 \end{equation}
Since there is no spatial variation in density, the flow is uniform along the pipe and we can take $q(x,t)=q(t)$. As water flows from one point to another, its pressure drops and the momentum equation can be used to determine the pressure loss, where we have
\begin{equation} \label{pde_water}
    \partial_t q(t)+A\partial_x p(x,t)+\rho Agsin\alpha+\frac{1}{\rho}\frac{\lambda}{2DA}q(t)|q(t)|=0,
\end{equation}
in which $\rho$ represents water density, $A$ is cross sectional area of the pipe, $g$ is gravitational acceleration, $\alpha$ is elevation of pipe, $D$ is diameter of the pipe and $\lambda$ is Darcy-Weisbach friction factor. The pressure gradient drives the flow and it remains constant along the pipe length\cite{jansen2013global}, so 
\begin{equation*}
    \partial_xp(0,t)=\frac{p(L,t)-p(0,t)}{L}.
\end{equation*}
This means that \eqref{pde_water} can be written as
\begin{equation*}
    \partial_t q(t)+A\frac{p(L,t)-p(0,t)}{L}+\rho Agsin\alpha+\frac{1}{\rho}\frac{\lambda}{2DA}q(t)|q(t)|=0
\end{equation*}
Using $S=\frac{L}{A}$, $H=-L\rho gsin\alpha$ and $g(q(t))=\frac{1}{\rho}\frac{\lambda Sq(t)|q(t)|}{2DA}$, we have the following pipeline equation
\begin{equation}\label{water pipe eqn}
    Sq^{'}(t)+p(L,t)-p(0,t)+g(q(t))=H
\end{equation} 

In case of a network of connected pipelines, one can use graph theory to model the network. Here we have two types of nodes,  pressure node $N_p$ and demand node $N_q$. To represent network mapping, the incidence matrix $A_G$ \cite{huck2014topology} is defined as:

\begin{equation*}
A_G=
    \begin{cases}
        \hspace{0.2cm} 1 & \text{if edge is incident at node} \\
       -1 & \text{if edge is leaving node } \\
        \hspace{0.2cm} 0 & \text{else }
    \end{cases}
\end{equation*}
We can split this incidence matrix into two sub-matrices representing only pressure node $A_G^p$ and demand node $A_G^q$ such that:
\begin{equation*}
    A_G \equiv A_G^q A_G^p
\end{equation*}
At node level, we have to make sure that the flow in and flow out for each node is equal to the demand flow. Mathematically we can describe it using mass balance:
 \begin{equation}\label{water node eqn}
     \sum_{I_{in}}q_i(t)-\sum_{I_{out}}q_j(t)=q_d(t)
 \end{equation}
 This set of equations define algebraic constraints for our network. All the information regarding water flow and demand at certain node is represented through the above equation. Combining these with \eqref{water pipe eqn}, we can represent the complete network model as system of differential algebraic equations (DAEs) system
\begin{equation} \label{odesys_water}
    \begin{aligned}
        S\frac{\partial q}{\partial t}+A_G^qp(t)+g(q(t))&=H-A_G^pp_s(t)\\
        (A_G^q)^Tq(t)&=q_s(t)
    \end{aligned}   
\end{equation}
where,
\[
S=\begin{bmatrix}
\frac{L_1}{A_1}   &         &           &   \\
                  &\frac{L_2}{A_2}      & \\
                  &         & \ddots    & \\
                  &         &           & \frac{L_i}{A_i}  \\
\end{bmatrix},~ q=\begin{bmatrix}
q_1  \\
q_2  \\
\vdots\\
q_i
\end{bmatrix}, ~p=\begin{bmatrix}
p_1  \\
p_2  \\
\vdots\\
p_i
\end{bmatrix}, ~g(q(t))=\begin{bmatrix}
\frac{L_1}{\rho}\frac{\lambda}{2DA^2_1}q_1(t)|q_1(t)|  \\
\frac{L_2}{\rho}\frac{\lambda}{2DA^2_2}q_2(t)|q_2(t)|  \\
\vdots\\
\frac{L_i}{\rho}\frac{\lambda}{2DA^2_i}q_i(t)|q_i(t)|
\end{bmatrix}
\]

\[
H=\begin{bmatrix}
-L_1\rho gsin\alpha_1\\
-L_2\rho gsin\alpha_2\\
\vdots\\
-L_i\rho gsin\alpha_i
\end{bmatrix}, \quad q_s(t)=\begin{bmatrix}
q_{s1}  \\
q_{s2}\\
\vdots\\
q_{si}
\end{bmatrix}, \quad p_s(t)=\begin{bmatrix}
p_{s1}  \\
p_{s2}\\
\vdots\\
p_{si}
\end{bmatrix}
\]
Clearly we can write \eqref{odesys_water} in matrix form as
\begin{equation}\label{power ode matrix}
        \begin{bmatrix} S & 0 \\0 & 0 \end{bmatrix}	
        \begin{bmatrix}  \partial_t q \\\partial_t p \end{bmatrix}=
        \begin{bmatrix} 0 & -A_G^p \\ (A_G^q)^T& 0 \end{bmatrix}
        \begin{bmatrix} q   \\ p \end{bmatrix}+
        \begin{bmatrix} -A_G^pp_s(t)   \\ -q_s(t) \end{bmatrix}+
        \begin{bmatrix} H-g(q(t))   \\ 0 \end{bmatrix}
\end{equation}
which is in the unified DAE form for water distribution networks
\begin{equation}\label{water ssp}
    E_w\dot{x}_w(t)=A_wx_w(t)+B_wu_w(t)+G_w
\end{equation}

\subsection{Power Distribution Network}
Electricity is transmitted from generation station to load centre and from load centre to consumers through transmission lines. For transmission line modeling, the distributed element model is used. This will help us to consider the spatial variance in electrical components like resistance, capacitance and inductance. Distributed element model assumes that values of these components spread continuously throughout circuit. Transmission line is divided into small section and each section have its own resistance, capacitance and inductance. Mathematically we can represent it using Telegrapher's equation \cite{bamigbola2014mathematical} in time domain.
\begin{equation}\label{telegraphers equation}
\begin{aligned}
    L\partial_t I+\partial_x V+RI&=0, \\
    C\partial_t V+\partial_x I+GV&=0,
\end{aligned}
\end{equation}
where $I$ is the distributed current, $V$ is the voltage across the unit, $R$ is the resistance, $L$ inductance, $C$ capacitance and $G$ conductance of the distributed elements. We used the finite difference method (FDM) for space discretization where standard forward difference results in
\begin{equation*}
    \frac{\partial V}{\partial x}=\frac{V_{i+1}-{V_i}}{\delta x},  \hspace{1cm}
    \frac{\partial I}{\partial x}=\frac{I_{i}-{I_{i-1}}}{\delta x},
\end{equation*}
which can be used to represent \eqref{telegraphers equation} as
\begin{equation}
     L_i\frac{\partial I_i}{\partial t}=-( -\frac{V_i}{\delta x}+\frac{V_{i+1}}{\delta x})-R_iI_i
\end{equation}
\begin{equation*}
    C_i\frac{\partial V_i}{\partial t}= -(\frac{I_i}{\delta x}-\frac{I_{i-1}}{\delta x})-G_iV_i=0
\end{equation*}
The discretized ODE model can be written in matrix form as
\begin{equation}\label{power ode matrix}
        \begin{bmatrix} L &  \\ & C \end{bmatrix}	
        \begin{bmatrix}  \partial_t I \\\partial_t V \end{bmatrix}=
        \begin{bmatrix} -R & -D_x \\ -D_x^T& -G \end{bmatrix}
        \begin{bmatrix} I \\ V \end{bmatrix}       
\end{equation}
In case of multiple transmission lines connected as a power distribution network, algebraic constraints for power network can be defined by using non-linear power flow equation \cite{gross2016solvability}. These equations make sure that the power injected into bus is equal to power flowing in network by following Kirchhoff current law (KCL). Our power network contains synchronous generator and two types of buses: generator buses in which case $P_i=p_{load}+p_{generator}$ and load buses where $P_i=p_{load}$. Total power in the $i-th$ bus becomes
\begin{equation*}\label{power generator}
    p_{generator}=\frac{E_i^{'}V_i}{X_i^{'}}[sin(\alpha_i-\theta_i)]
\end{equation*}
Now, taking non-linear power flow equation which will help us to study natural power flow in network:
\begin{equation}
    P_i=\sum_{j=1}^{n+m}|V_i||V_j|[G_{ij}cos(\theta_i-\theta_j)+B_{ij}sin(\theta_i-\theta_j)]
\end{equation}
Then, algebraic constraint for our system
\begin{equation}
    \begin{aligned}
        p_{load}+p_{generator}&=P_i  \hspace{1cm}\text{(for generator bus where $p_{generator}$ from \eqref{power generator}) }\\
        p_{load}&=P_i \hspace{1cm}\text{(for load bus) }
    \end{aligned}
\end{equation}
Writing above equation in form:
\begin{equation}\label{power constraint}
    p-G(V,\theta)=0
\end{equation}
Where, $G(V,\theta)$ represents non linearity in system.
Equations \eqref{power ode matrix} and \eqref{power constraint} combine to form the differential algebraic equation system for power distribution network i.e.:
\begin{equation}\
        \begin{bmatrix} L &  \\ & C \\ & & 0 \end{bmatrix}	
        \begin{bmatrix}  \partial_t I \\\partial_t V \\ \partial_t \theta \end{bmatrix}=
        \begin{bmatrix} -R & -D_x &0 \\ -D_x^T& -G &0 \\0&0&0 \end{bmatrix}
        \begin{bmatrix} I   \\ V \\\theta \end{bmatrix}+
        \begin{bmatrix} 0   \\ 0 \\ p \end{bmatrix}+\begin{bmatrix} 0   \\ 0 \\ -G(V,\theta) \end{bmatrix}      
\end{equation}
where,
\[ \footnotesize 
L=\begin{bmatrix}
  L_1 &         &           &     \\
                  &L_2      &    \\
                  &         & \ddots        &     \\
                  &         &           &  L_i 
\end{bmatrix}
;~C=\begin{bmatrix}
  C_1 &         &           &    \\
                  &C_2      &     \\
                  &         & \ddots       &  \\
                  &         &           &  C_i  \\
\end{bmatrix}
;~R=\begin{bmatrix}
  R_1 &         &           &     \\
                  &R_2      &     \\
                  &         & \ddots         &    \\
                  &         &           & R_i  \\
\end{bmatrix}
\]
\[ \footnotesize 
G=\begin{bmatrix}
  G_1 &         &           &     \\
                  &G_2      &    \\
                  &         & \ddots        &     \\
                  &         &           & G_i  \\
\end{bmatrix}, \quad 
D_x=\frac{1}{\Delta x}\begin{bmatrix}
  -1 &  1       &           & &    \\
     &-1        &1          & &   \\
     &          &\ddots        &\ddots&    \\
     &         &           & -1&1  
\end{bmatrix}. 
\]
State space representation of the power distribution network is
\begin{equation}\label{power ssp}
    E_p\dot{x}_p(t)=A_px_p(t)+B_pu_p(t)+G_p
\end{equation}

\section{Model Reduction of the Unified Network}\label{s3}
In this section we show how to construct a reduced system of the form \eqref{rom} for a given DAE system as in \eqref{fom}.
The concept of projection can be used to construct the reduced model, where the state vector $x(t)$ is projected to an $r$-dimensional subspace spanned by columns of $V\in \mathbb{R}^{n\times r}$, and it is assumed that $x(t)\approx Vx_r(t)$ for some $x_r \in \mathbb{R}^{r}$. In case of Petrov-Galerkin projection, it is assumed that the residual error associated with the approximation is orthogonal to another $r$ dimensional subspace spanned by the columns of $W\in \mathbb{R}^{n\times r}$ (Galerkin projection if $W=V$).
\begin{equation*}\label{reducedssp}
    W^TR(x)=W^T\Big(EV\dot{x}_r(t)-AVx_r(t)-Bu(t)-G(x_r,u)\Big)=0, \quad \mbox{(Petrov-Galerkin)}
\end{equation*}
This mean that by using Petrov-Galerkin Condition, we can construct the reduced systems in terms of basis matrices $V$ and $W$. The important question is how to choose $V$ and $W$. In case of linear systems, where the nonlinear part $G(x,u)=0$, these basis matrices are linked to the transfer function of the linear system. The transfer function of the system in \eqref{fom} with $G(x,u)=0$ and the corresponding reduced system's transfer function can be written as
\begin{equation*}
    H(s)=C(sE-A)^{-1}B, \quad
    H_r(s)=C_r(sE_r-A_r)^{-1}B_r.
\end{equation*}
It is shown in the literature \cite{} that in case of single input single output (SISO) linear systems with nonsingular $E$, the specific choice of basis matrices $V$ and $W$      
\begin{equation*}
    \begin{aligned}
        V&=[(\sigma_1E-A)^{-1}B,...,(\sigma_rE-A)^{-1}B],\\
        W&=[(\sigma_1E-A)^{-T}C^T,...,(\sigma_rE-A)^{-T}C^T],
    \end{aligned}
\end{equation*}
ensure that the associated reduced system satisfy the following interpolation conditions
\begin{equation*}
    H(\sigma_i)= H_r(\sigma_i), \quad i = 1,\ldots, r
\end{equation*}
This means that the reduced system's transfer function exactly matches the original transfer function at some fixed interpolation points $\sigma_i, ~i=1,\ldots,r$. The question is how to select these interpolation points. In \cite{gugercin2013model}, a Newton iteration-based framework has been proposed that iteratively updates the interpolation points until they converge to the mirror images of the eigenvalues of the reduced system. On convergence, the reduced system is optimal among all r-dimensional reduced systems in terms of $\mathcal{H}_2$ norm. That is 
\begin{equation*}
    H_r = \arg \min_{dim(\tilde{H})=r} \|H(s)-\tilde{H}(s)\|_{\mathcal{H}_2}
\end{equation*}
This algorithm is called the iterative rational Krylov algorithm (IRKA) and it is well used in the literature for model order reduction of SISO linear systems. The approach has been extended to multi-input multi-output (MIMO) systems, where the interpolation conditions are satisfied in some tangent direction, $b_i\in \mathbb{R}^m$ and $c_i\in \mathbb{R}^p$ for $i = 1,\ldots, r$. That is
\begin{equation*}
    H(\sigma_i)b_i= H_r(\sigma_i)b_i,~\mbox{ or }~ c_i^TH(\sigma_i)= c_i^TH_r(\sigma_i),~\mbox{ or }~  c_i^TH(\sigma_i)b_i= c_i^TH_r(\sigma_i)b_i.
\end{equation*}
The IRKA type framework for tangential interpolation has been developed in the literature \cite{druskin2014adaptive} and is called tangential IRKA or tIRKA. The algorithm can reduce MIMO systems and in case of SISO systems, it automatically reduces to the standard IRKA. Extension of these methods to the reduction of linear descriptor systems or system of DAEs is also possible if the associated matrix pencil $(E,A)$ is stable. The issue however is that the reduced system is strictly proper ODE system for a given DAE system. This means that at large frequencies the error system will theoretically goes to infinity because the improper part of the original DAE can reach infinity. To overcome this issue, \cite{} utilizes the concept of spectral projectors that explicitly identify the strictly proper and improper polynomial part of the system and the IRKA type method reduce the strictly proper part while retaining the polynomial part without reduction. For some structured descriptor systems, the polynomial part is constant and can be retained in the reduced system by the $D_r$-term. It is shown in \cite{} that for some matrices $F \in \mathbb{R}^{n \times m}$ and $\bar{F}\in \mathbb{R}^{n \times m}$ satisfying $F^TV=\begin{bmatrix} b_1, \cdots, b_r \end{bmatrix}$ and $W^T\bar{F}=\begin{bmatrix} c_1, \cdots, c_r \end{bmatrix}^T$, then for a reduced system with transfer function 
$$
    \hat{H}(s)=\hat{C}(sE_r-\hat{A})^{-1}\hat{B}+D_r,
$$
in which $\hat{A}= A_r-\bar{F}D_rF^T$, $\hat{B}=B-\bar{F}D_r$ and $\hat{C}= C_r-D_rF^T$, the following interpolation condition holds:
\begin{equation*}
    H(\sigma_i)b_i= \hat{H}(\sigma_i)b_i,~\mbox{ and }~ c_i^TH(\sigma_i)= c_i^T\hat{H}(\sigma_i),~\mbox{ and }~  c_i^TH(\sigma_i)b_i= c_i^T\hat{H}(\sigma_i)b_i.
\end{equation*}
This linear model reduction framework can be utilized for the reduction of the nonlinear DAE system by constructing $V$ and $W$ from the linear part of the nonlinear DAE system \eqref{fom} and projecting it with these linear basis matrices. In addition, the polynomial part should be retained in the non-zero $D_r$ term. The complete framework is summarized in the following algorithm:
\begin{algorithm}
\caption{tangential IRKA with nonzero $D_r$}\label{alg:cap}
\begin{algorithmic}
\Require $E,~A,~B,~C,~r,~S_r=\sigma_i|_1^r,~\mathcal{B}=b_i|_1^r,~\mathcal{C}=c_i|_1^r$
\State Define 'tol' and 'max-iter'
\While{$err < \mbox{tol}$ and iter $<$ max-iter}
\State $   V\gets \begin{bmatrix}(\sigma_1E-A)^{-1}Bb_1&\cdots &(\sigma_rE-A)^{-1}Bb_r\end{bmatrix}$
\State $   W\gets \begin{bmatrix}(\sigma_1E-A)^{-T}C^Tc_1&\cdots &(\sigma_rE-A)^{-T}C^Tc_r\end{bmatrix}$
\State Compute $A_r = W^TAV$
\State $error = \|sort(\lambda(-A_r)) - sort(S_r)\|$
\State Update $S_r \gets \lambda(-A_r)$ and tangent directions
\EndWhile
\State Solve $F^TV=\mathcal{B}$ and $W^T\bar{F}=\mathcal{C}^T$ for $F$ and $\bar{F}$
\State $E_r\gets W^TEV$,~ $\hat{A}\gets W^TAV+\bar{F}D_rF^T$,~ $\hat{B}\gets W^TB-\bar{F}D_r$,~\\ $\hat{C}\gets CV-D_rF^T$ and $G_r(x_r,u)=W^TG(Vx_r,u)$
\end{algorithmic}
\end{algorithm}

\section{Numerical Results}
We present two different discretization methods, the finite volume method (FVM) and the finite difference method (FDM), to model the behavior of a simple gas distribution network and show how the reduction framework given in Algorithm~1 can be used to efficiently simulate the model. 
The comparison of the discretization method is done in terms of execution time, pressure and mass flow responses.

\subsection{Single Pipeline Connections Network}
\begin{figure}[H]
    \centering
    \includegraphics[width=0.75\linewidth]{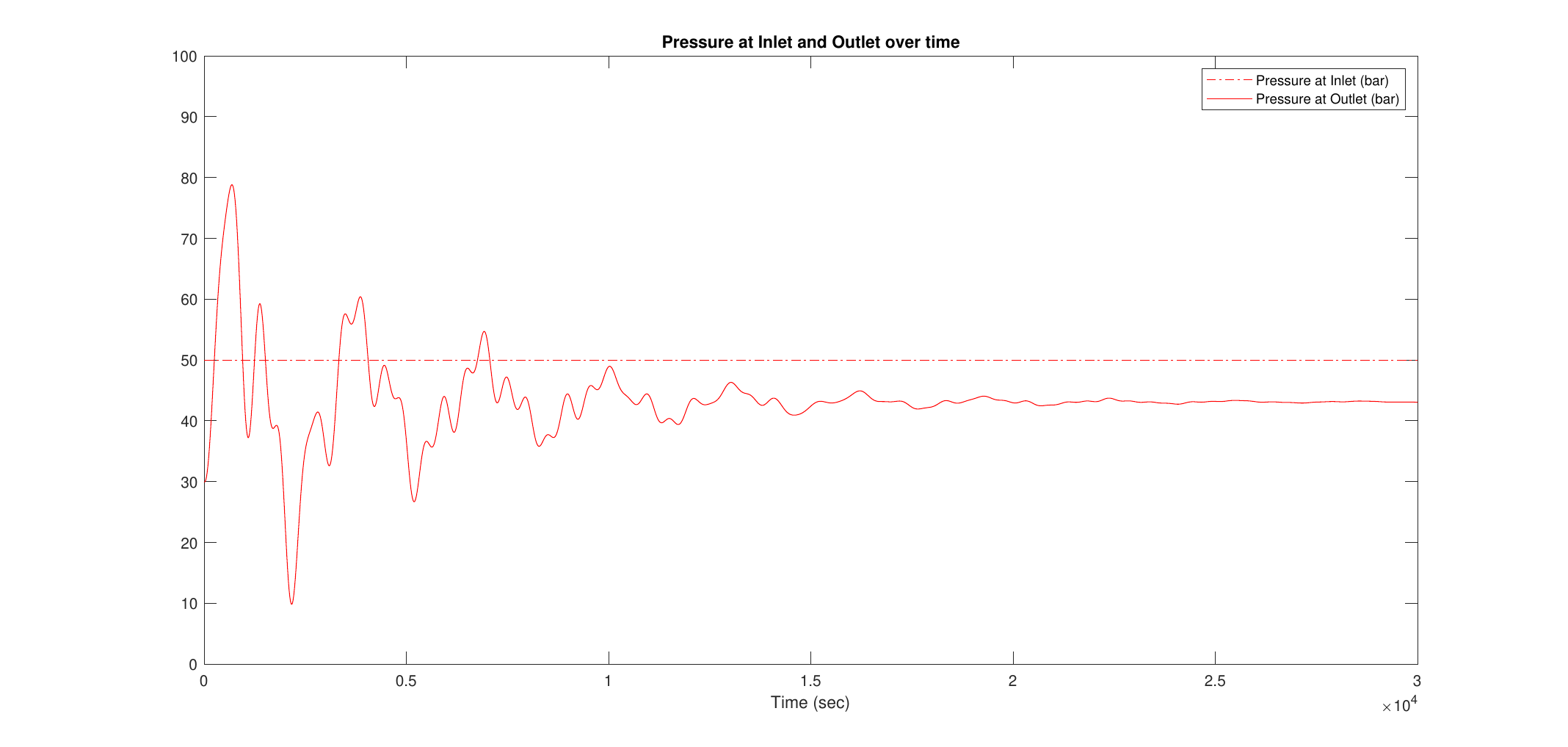}
    \caption{FVM Pressure Profile}
    \label{fig:fvmPressure}
\end{figure}

\begin{figure}[H]
    \centering
    \includegraphics[width=0.75\linewidth]{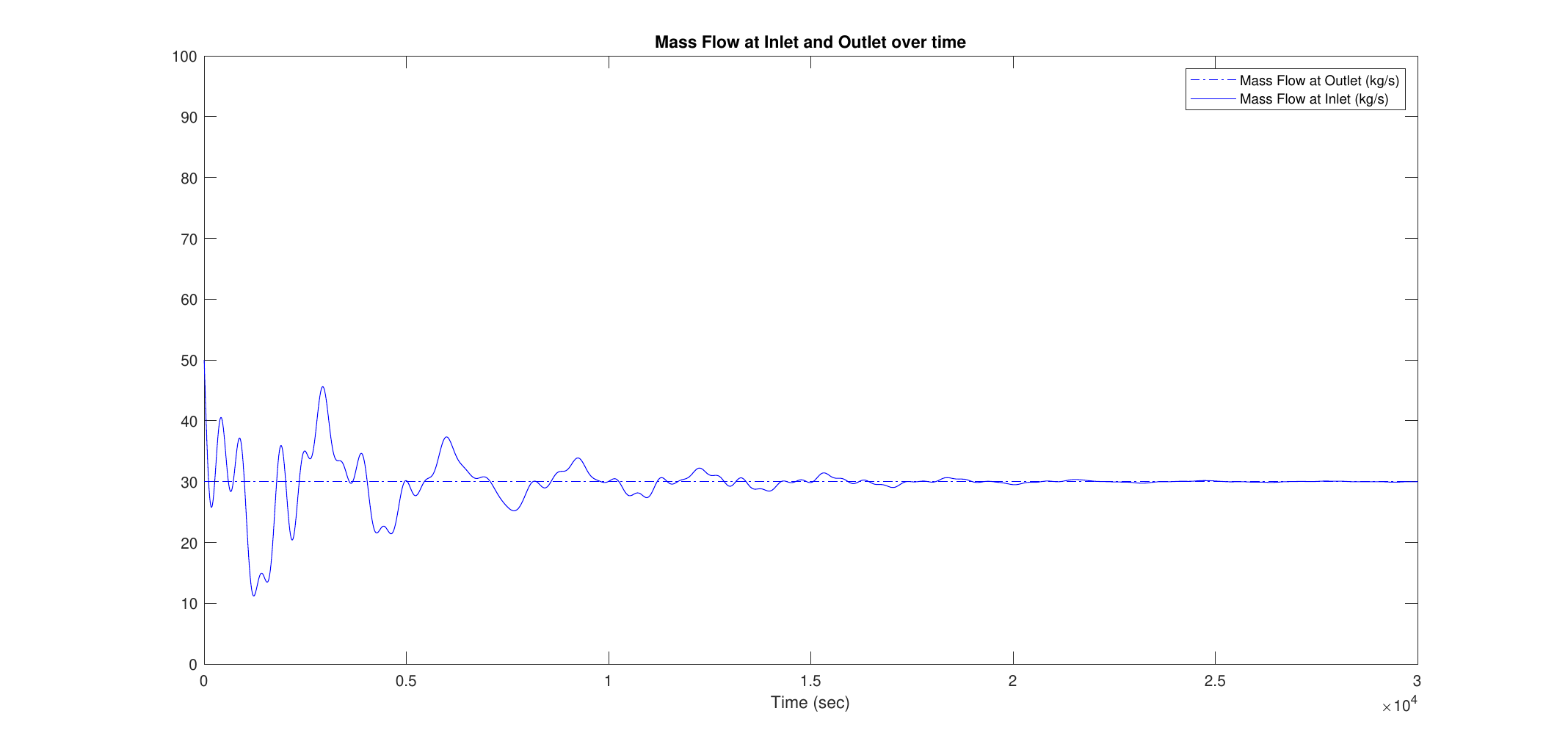}
    \caption{FVM Mass Flow Profile}
    \label{fig:fvmmassflow}
\end{figure}

The first example represents a benchmark model for analysis of gas distribution network \cite{qiu2018efficient} where a single pipeline connection is allowed at each node. The network is defined in terms of the physical specifications given in Table~\ref{tab:sys_specs}. The mathematical representation of the network is made using the finite-volume method as discussed in Section~2, resulting in a discretized system of the form \eqref{mat_vec_form}. The Dirichlet boundary conditions are defined such that on the inlet side the supplied pressure is $p_s=50$ and on the outlet side the demand mass flow is $q_d=30$. With these boundary conditions, the ODE system is solved through the implicit Euler method, where if the mass matrix is represented by $E_v$ and time step as $\tau$, we can iteratively converge to the required solution as
\begin{equation}\label{implicitEuler}
    x_v^k = (E_v-\tau A_v)^{-1}(E_vx_v^{k-1}+\tau B_vu - \tau G_v).
\end{equation}
The matrix pencil is $(E_v-\tau A_v)$ is assumed to be invertible and the iterations $k$ are bounded by $\mbox{max-iter}=1000$. The matrices with subscript $v$ shows that they are constructed through the FVM method. The pressure at outlet of the gas network and the mass flow at inlet of the network obtained through the solution of the ODE system is shown in Figure \ref{fig:fvmPressure} and \ref{fig:fvmmassflow}. The results show that after some variation in the pressure, it settles at a value lower than the supplied pressure setting the supplied mass flow equal to the required demand mass flow at the outlet of the network. 
\begin{table}[H]
    \centering
    \begin{tabular}{|p{6cm}|p{3cm}|}
    \hline
    \textbf{Parameters}& \textbf{Measurements}\\
    \hline\hline
         Pipe length     &1000m \\
         Diameter of pipe&1.0 m \\
         Area of pipe    &$0.7854 m^2$ \\
         Supply pressure, $p_s$ &50 bar \\
         Demand Flow, $q_d$     &30 kg/s \\
         Mesh Size       &100m \\
    \hline
    \end{tabular}
    \caption{System specifications}
    \label{tab:sys_specs}
\end{table}

\begin{figure}[H]
    \centering
    \includegraphics[width=0.75\linewidth]{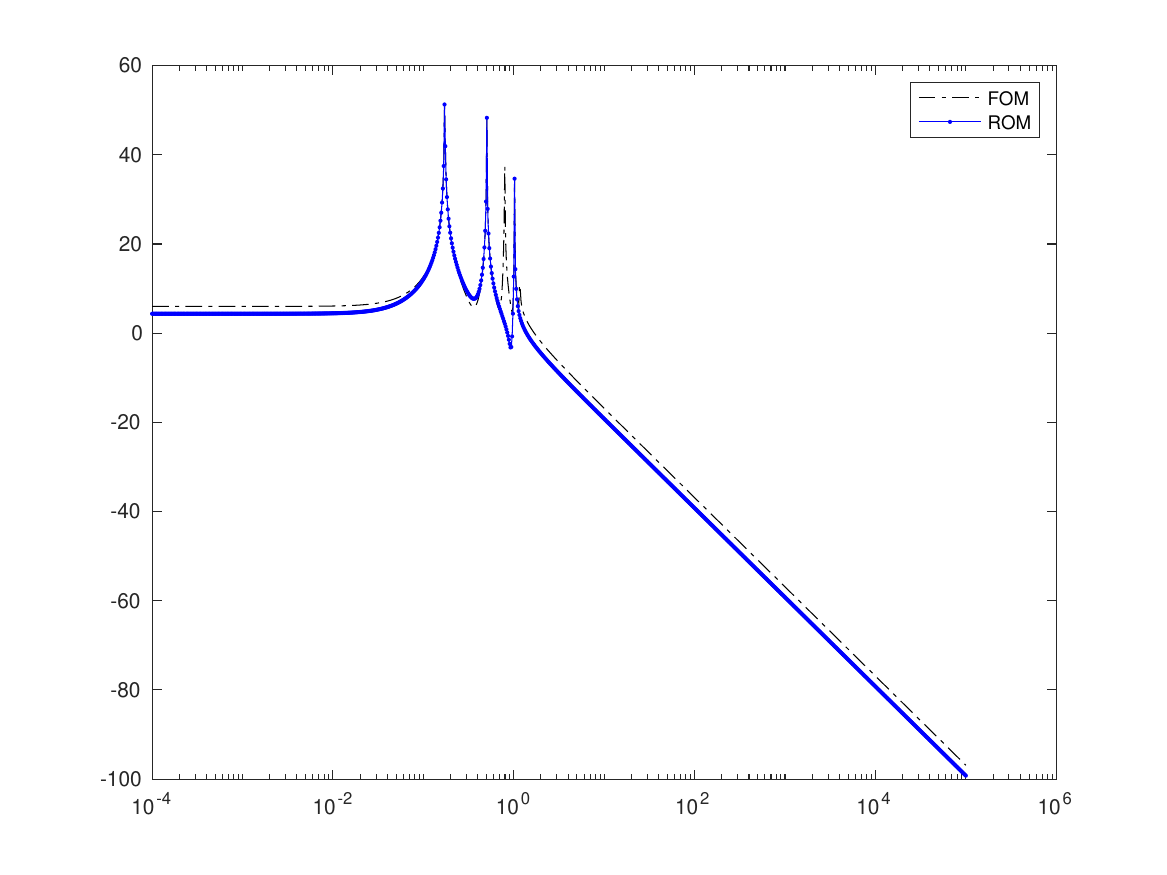}
    \caption{MOR of the linear part}
    \label{fig:mor}
\end{figure}

The model order reduction of the linear part is done as discussed in Section~3. The maximum singular value of the transfer matrix for original and reduced system are shown in Figure \ref{fig:mor} for $r=6$ and $2n-2=18$.

Similarly, we can perform the finite difference method (FDM) on the governing equations to obtain a nonlinear matrix equation similar to \eqref{mat_vec_form} but with $E$ an identity matrix now. Keeping all the parameters as defined before, the implicit Euler method can be written as :
\begin{equation}
    x_d^k = (I_d-\tau A_d)^{-1}(x_d^{k-1}+\tau B_du-\tau G_d)
\end{equation}
\begin{figure}[H]
    \centering
    \includegraphics[width=0.75\linewidth]{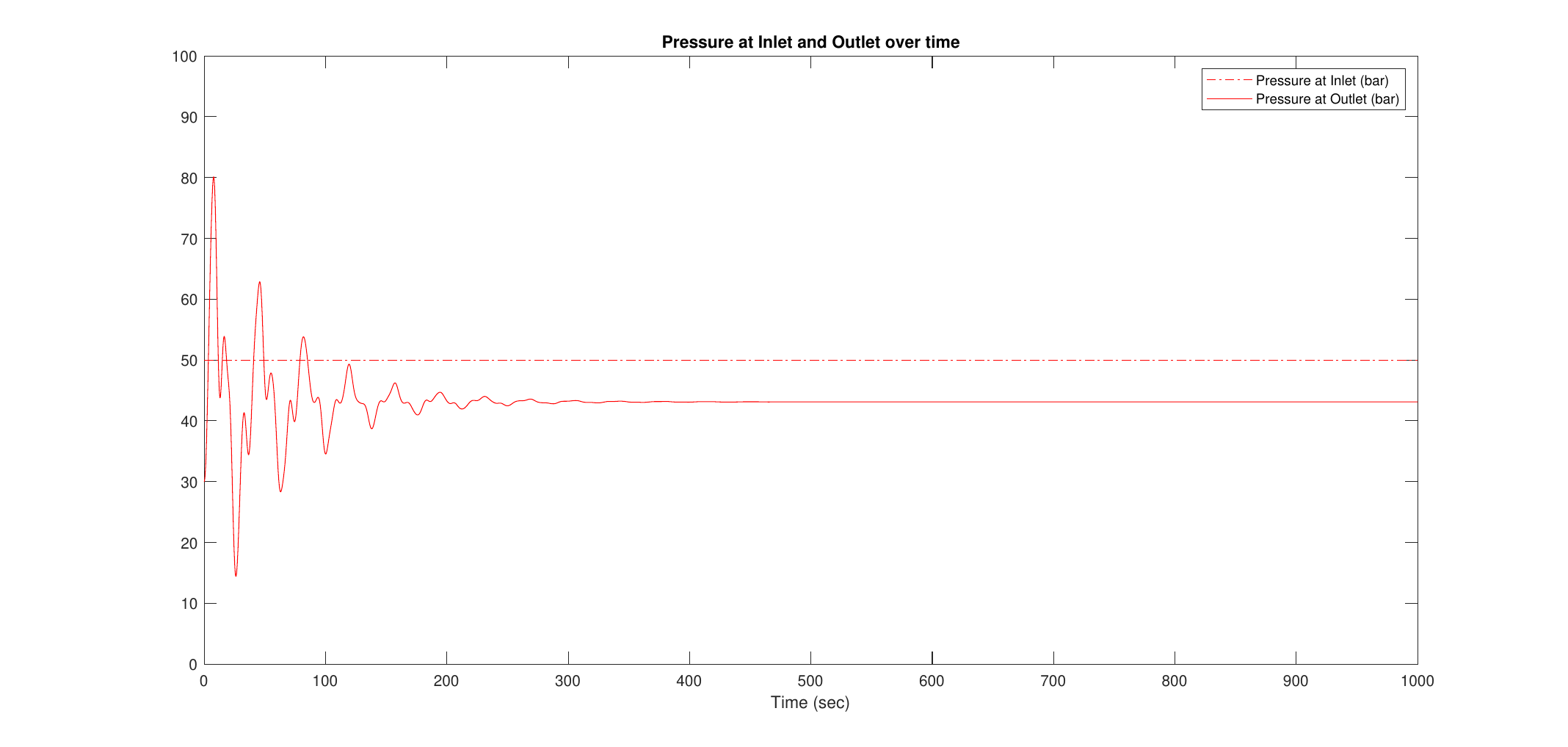}
    \caption{FDM Pressure Profile}
    \label{fig:fdmpressure}
\end{figure}
\begin{figure}[H]
    \centering
    \includegraphics[width=0.75\linewidth]{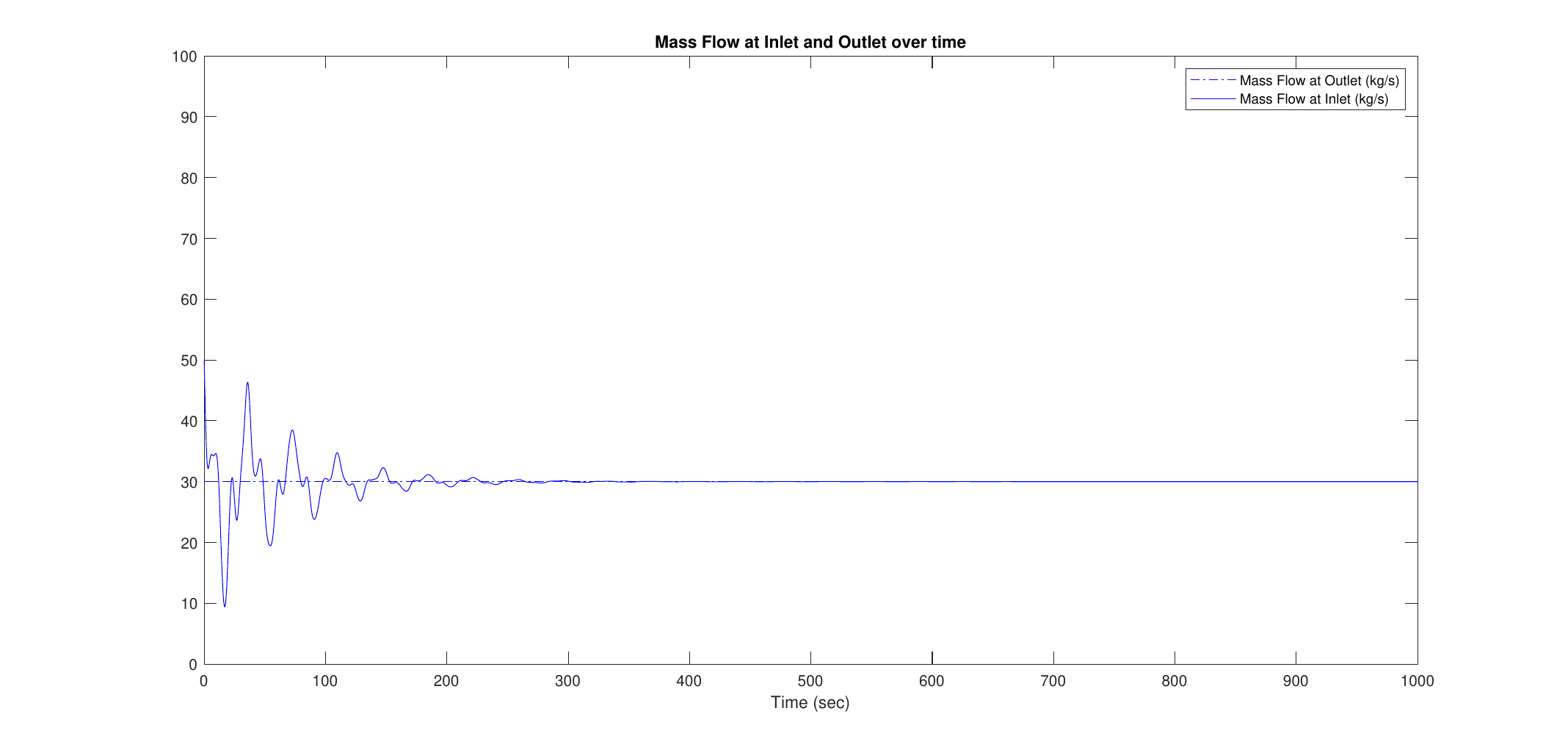}
    \caption{FDM Mass Flow Profile}
    \label{fig:fdmmassflow}
\end{figure}

The results show that the modeling approach captures the fundamental characteristics of the network such as time-dependent changes in pressure and mass flow due to frictional losses. In addition, mass conservation can be observed as inlet and outlet flows tend to balance over time.  In Figure \ref{fig:fdmpressure} pressure profile at outlet starts with strong oscillations which corresponds to wave propogation but settles due to frictional dissipation. Final equilibrium pressure emerges as system adjusts to match inlet pressure and outlet demand as per requirement of conservation laws. The mass flow profile in Figure \ref{fig:fdmmassflow}  at inlet shows oscillations as pipeline's initial attempt to balance the pressure differences. System stabilizes and steady- state balance is achieved as inlet mass flow equals 30kg/s.

All the simulations are performed on (mention system specs). Table \ref{tab:exec_time} below gives execution time for different time steps.
\begin{table}[H]
\centering
    \begin{tabular}{|p{2cm}|p{4cm}|p{4cm}|}
        \hline
        Time Step (sec) & Execution Time FVM (sec) &Execution Time FDM (sec)\\
        \hline
         1.00  &0.866  &0.948\\
         0.50  &1.297  &1.444\\
         0.25  &2.284  &2.304\\
         0.10  &5.081  &5.155\\
         0.05  &10.017 &10.049\\
         0.01  &48.995 &48.749\\
         0.001 &491.278&497.246\\
         \hline
    \end{tabular}
    \caption{Execution time at different time steps}
    \label{tab:exec_time}
\end{table}

\section{Conclusion}
In conclusion, this research work highlights the importance of discretization methods and model order reduction techniques in the simulation of complex energy distribution networks. We demonstrated that the energy networks such as gas, water, and power distribution network can be represented by a unified system of nonlinear differential algebraic equations (DAE) and it can be reduced through a new tangential IRKA type model reduction framework that retains the improper part of the original system allowing bounded value for the error system especially at large frequencies. The reduced system closely approximates the behavior of the original unified mathematical model while significantly decreasing the online computational costs.

\paragraph{Author Contributions.} The Authors (Saleha Kiran, Farhan Hussain and Mian Ilyas Ahmad) have contributed as follows: Saleha Kiran: MATLAB implementation and testing, Farhan Hussain: writing – original draft preparation and review, and Mian Ilyas Ahmad: writing – review and editing. All authors have read and approved the published version of the manuscript.


\paragraph{Conflicts of Interest.} The authors declare no conflicts of interest.


\section*{References}
\bibliographystyle{elsarticle-num.bst} 
\bibliography{QBDAE_2s} 

\end{document}